\RequirePackage[leqno]{amsmath}
\documentclass[leqno]{amsart}
\usepackage{amsmath}
\makeatletter
\usepackage[utf8]{inputenc}     % Accept different input encodings
\usepackage{multicol}
\usepackage{enumitem}
\usepackage{color, graphics}
\usepackage[demo]{graphicx}
\usepackage{efbox,graphicx}
\efboxsetup{linecolor=gray,linewidth=0.8pt}
\usepackage{floatrow}
\usepackage{pgfplots}
\usepackage{graphicx, multicol, latexsym, amsmath, amssymb}
\usepackage{pgfplots}
\usepackage{tikz}
\usepackage{amsfonts}
\usepackage{enumitem}
\setlist[itemize]{topsep=0pt,after=\vspace{1.5\baselineskip}}
 \usepackage[colorlinks=true]{hyperref}
\usepackage{empheq}
\usepackage[T1]{fontenc}
\usepackage{tabularx}
\usepackage[leftcaption]{sidecap}
\sidecaptionvpos{figure}{m}
\usepackage{tikz}
\usepackage{subfigure}
\usepackage[margin=1in]{geometry}
\usetikzlibrary{backgrounds,automata}

\newcount\rc@count
\let\rc@clearconstantlist\empty
\newcommand\rc@clearconstant[1]{\global\expandafter\let\csname rc@const@#1\endcsname\undefined}
\newcommand\resetconstants[1]{%
    \def\rc@constname{#1}% Set the new base name of the constants to the argument
    \global\rc@count=1\relax % Reset the constant counter to 1
    \bgroup 
        \let\\\rc@clearconstant % map over the list of constants that have been defined, clearing each of them.
        \rc@clearconstantlist
        \global\let\rc@clearconstantlist\empty % Globally empty the list of constants.
    \egroup
}
\newcommand\const[1]{%
    \@ifundefined{rc@const@#1}{%
        \expandafter\xdef\csname rc@const@#1\endcsname{%
           \noexpand\rc@useconst{\rc@constname}{\the\rc@count}%
        }%
        \g@addto@macro\rc@clearconstantlist{\\{\mathrm{#1}}}%
        \global\advance\rc@count1\relax
    }{}%
    \csname rc@const@#1\endcsname
}
\newcommand\rc@useconst[2]{{#1}\textsubscript{#2}}
\setlist[itemize]{noitemsep, topsep=0pt}

\def\R{\mathbb R} \def\N{\mathbb N}

\def\R{\mathbb R} \def\N{\mathbb N} 
\def\TM{T_{max}} 
\def
\@cite
#1#2{[#1\if@tempswa, #2\fi]}
\makeatother
\newtheorem{theorem}{Theorem}[section]

\newtheorem{lemma}[theorem]{Lemma}

\newtheorem{remark}{Remark}

\newcounter{cnstcnt}

\title[Nonlocal effects on a chemotaxis model]{Boundedness through nonlocal dampening effects in a fully parabolic chemotaxis model with sub and superquadratic growth}
\author[Yutaro Chiyo, Fatma Gamze D\"{u}zg\"{u}n, Silvia Frassu and Giuseppe Viglialoro]{}
\keywords{Chemotaxis, Global existence, Nonlocal growth terms, Boundedness. \\
\textit{$^*$Corresponding author}: silvia.frassu@unica.it}

\begin{document}

\maketitle
%\tableofcontents

\centerline{\scshape{\scshape{Yutaro Chiyo$^{\natural}$, Fatma Gamze D\"{u}zg\"{u}n$^{\flat}$, Silvia Frassu$^{\sharp,*}$ \and Giuseppe Viglialoro$^{\sharp}$}}}
\medskip
{
\medskip
\centerline{$^{\natural}$Department of Mathematics} 
\centerline{Tokyo University of Science}
\centerline{1-3, Kagurazaka, Shinjuku-ku, Tokyo 162-8601 (Japan)}
\medskip
}
{
\medskip
\centerline{$^{\flat}$Department of Mathematics} 
\centerline{Hacettepe University}
\centerline{06800, Beytepe, Ankara (Turkey)}
\medskip
}
{
\medskip
\centerline{$^\sharp$Dipartimento di Matematica e Informatica}
\centerline{Universit\`{a} di Cagliari}
\centerline{Via Ospedale 72, 09124. Cagliari (Italy)}
\medskip
}
\bigskip
\begin{abstract}
This work deals with a chemotaxis model where an external source involving a sub and superquadratic growth effect contrasted by nonlocal dampening reaction influences the motion of a cell density attracted by a chemical signal. We study the mechanism of the two densities once their initial configurations are fixed in bounded impenetrable regions; in the specific, we establish that no gathering effect for the cells can appear in time provided that the dampening effect is strong enough. 

Mathematically, we are concerned with this problem 
\begin{equation}\label{problem_abstract}
\tag{$\Diamond$}
\begin{cases}
u_t=\Delta u-\chi\nabla \cdot (u\nabla v)+au^\alpha-bu^\alpha\int_\Omega u^\beta &{\rm in}\ \Omega \times (0, \TM),\\
\tau v_t=\Delta v-v+u &{\rm in}\ \Omega \times (0, \TM),\\
u_\nu=v_\nu=0 &{\rm on}\ \partial\Omega \times (0, \TM),\\
u(x, 0)=u_0(x)\geq 0, v(x,0)=v_0(x)\geq 0, &x \in \bar{\Omega},
\end{cases}
\end{equation}
for $\tau=1$, $n\in \N$, $\chi,a,b>0$ and $\alpha, \beta\geq 1$. Herein $u$ stands for the population density,  $v$ for the chemical signal and $\TM$ for the maximal time of existence of any nonnegative classical solution $(u,v)$ to system \eqref{problem_abstract}. We prove that despite any large-mass initial data $u_0$, whenever 
\begin{itemize}
\item (the subquadratic case) $1\leq \alpha <2 \quad \textrm{and} \quad \beta>\frac{n+4}{2}-\alpha,$
\item (the superquadratic case) $\beta>\frac{n}{2} \quad \textrm{and} \quad 2\leq \alpha < 1+ \frac{2\beta}{n},$
\end{itemize}
actually $\TM=\infty$ and $u$ and $v$ are uniformly bounded.

This paper is in line with the result in \cite{Bian2018}, where the same conclusion is established for the simplified parabolic-elliptic version of model \eqref{problem_abstract}, corresponding to $\tau=0$; more exactly, this work extends the study to the fully parabolic case \cite{Bian2018}.
\end{abstract}
\resetconstants{c}
\tableofcontents
\section{Introduction and motivations}\label{IntroSection}
\subsection{Basic description of the research}
In this paper we consider 
\begin{equation}\label{problem}
\begin{cases}
u_t=\Delta u-\chi\nabla \cdot (u\nabla v)+au^\alpha-bu^\alpha\int_\Omega u^\beta &{\rm in}\ \Omega \times (0, \TM),\\
v_t=\Delta v-v+u &{\rm in}\ \Omega \times (0, \TM),\\
u_\nu=v_\nu=0 &{\rm on}\ \partial\Omega \times (0, \TM),\\
u(x, 0)=u_0(x), v(x,0)=v_0(x) &x \in \bar{\Omega},
\end{cases}
\end{equation}
where $\Omega \subset \mathbb{R}^n$ ($n\in \N$) is a bounded domain with smooth boundary $\partial\Omega$ (briefly, ``bounded and smooth domain''); additionally, we fix  $\chi, a, b>0$, $\alpha, \beta \ge 1$ and sufficiently regular and nonnegative initial data $u_0(x), v_0(x)$. On the other hand, the subscript $\nu$ in $(\cdot)_\nu$ indicates the outward normal derivative on $\partial\Omega$ and $\TM$ is the maximal existence time up to 
which solutions to the system are defined.

If properly interpreted, this model idealizes a \emph{chemotaxis} phenomenon, a mechanism from mathematical biology describing the directed migration of a cell in response to a chemical signal; more exactly, the movement of an organism or entity (such as somatic cells, bacteria, and other
single-cell or multicellular organisms) is strongly influenced by the presence of a stimulus, and precisely the motion follows the direction of the gradient of the stimulus itself. 

It is well known that the land marking event of chemotaxis was first introduced by Keller and Segel in 1970s (\cite{Keller1970}, \cite{Keller-Segel1970}). More expressly, by indicating with  $u = u(x, t)$ a certain cell density at the position $x$ and at the time $t$, and with $v=v(x,t)$ the stimulus at the same position and time, the pioneering study reads as \eqref{problem} for the specific case $a=b=0$.  The partial differential equation modeling the motion of  $u$, i.e. 
\begin{equation}\label{problemOriginalKS}
u_t=\Delta u-\chi  \nabla \cdot (u \nabla v)\quad  \textrm{in}\quad \Omega \times (0,\TM),
\end{equation}
essentially describes how a chemotactical impact of the (chemo)sensitivity ($\chi$) provided by the chemical signal $v$  may break the natural diffusion  (associated to the Laplacian operator, $\Delta u$) of the cells. Indeed, the term $-\nabla \cdot (u \chi \nabla v)$ models the transport of $u$ in the direction  $\chi \nabla v$, the negative sign indicating the attractive effect that $v$ has  on the cells (higher for $\chi$ larger and for an increasing amount of $v$). As a consequence, when $v$ is produced by the same cells, and in such a scenario $v$ obeys 
\begin{equation}\label{SecondEqKS}
v_t=\Delta v-v+u\quad \textrm{in}\quad \Omega \times (0,\TM),
\end{equation}
the attractive impact may be so efficient as to lead the cell density to its chemotactic collapse (blow-up at finite time with appearance of $\delta$-formations in the region).
\subsection{An overview on the Keller--Segel system}
Mathematically, it was proved that solutions to the initial-boundary value problem associated to equations \eqref{problemOriginalKS} and \eqref{SecondEqKS}, may be globally bounded in time or may blow up at finite time; this depends on the mass (i.e., $\int_\Omega u_0(x)dx$) of the initial data, its specific configuration, and the value of the sensitivity $\chi$. More precisely, in one-dimensional settings, all solutions are uniformly bounded in time, whereas for $n\ge 3$ given any arbitrarily small mass $m=\int_\Omega u_0(x)dx>0$, it is possible to construct solutions  blowing-up at finite time. On the other hand, when $n=2$, the value $4 \pi$ separates the case where diffusion overcomes self-attraction (if $\chi m<4\pi$) from the opposite scenario where self-attraction dominates (if $\chi m>4\pi$); respectively, all solutions are global in time, and initial data producing assembling processes at finite time can be detected. A detailed discussion on such analyses can be found in \cite{HerreroVelazquez,Nagai,Nagai-Radial,win_aggregation_vs}, which are undoubtedly classical results in this context. 
\subsection{An overview on the Keller--Segel system with logistics}
If the evolution of $u$ in equation \eqref{problemOriginalKS} is also influenced by the presence of logistic terms behaving as $au-bu^{\beta}$, for $\beta>1$, mathematical intuition suggests that superlinear damping effects should benefit the boundedness of solutions (this, for instance, occurs for ordinary differential equations of the type $u'=au-bu^{\beta}$). Actually,  the prevention of $\delta$-formations in the sense of finite-time blow-up for 
\begin{equation}\label{KS}
u_t=\Delta u-\chi  \nabla \cdot (u \nabla v)+au-bu^\beta\quad  \textrm{in}\quad \Omega \times (0,\TM),
\end{equation}
when coupled with some equation implying the segregation of $v$ with $u$ (for instance \eqref{SecondEqKS}), has been established only for large values of $b$  (if $\beta=2$, see \cite{TelloWinkParEl}, \cite{W0}), whereas for some value of $\beta$ near $1$ a blow-up scenario was detected, first for dimension $5$ or higher \cite{WinDespiteLogistic}, (see also \cite{FuestCriticalNoDEA} for an improvement of \cite{WinDespiteLogistic}), but later also for $n\geq 3$, in  \cite{Winkler_ZAMP-FiniteTimeLowDimension}.
\subsection{An overview on the Keller--Segel system with nonlocal sources}
As anticipated, in this research we are interested in understanding how the introduction of external 
growth factors of logistic type defined in terms of the total mass of the some power of the population, and hence idealized by nonlocal external sources,  may avoid blow-up mechanisms, exactly as  logistics. 
To be precise, likewise to classical logistic effects, impacts behaving as
\begin{equation}\label{NonLocalSources}
au^\alpha-bu^\alpha\int_\Omega u^\beta \quad a,b>0 \textrm{ and } \alpha,\beta\geq 1,
\end{equation}
model a competition between a birth contribution, favoring instabilities of the species (especially for large values of $a$), and a death one opportunely contrasting this instability (especially for large values of $b$).

In this context, some questions naturally arise. 
\begin{enumerate}[label=\textbf{$\mathcal{Q}$:},ref=$\mathcal{Q}$]
%\item a
\item \label{Question} Can one expect that in a biological mechanism governed by the equation
\begin{equation}\label{KS-NonLocal}
u_t=\Delta u-\chi  \nabla \cdot (u \nabla v)+au^\alpha-bu^\alpha\int_\Omega u^\beta \quad  \textrm{in}\quad \Omega \times (0,\TM),
\end{equation}
the external dampening source suffices to enforce boundedness of solutions, even for any large initial distribution $u_0$, arbitrarily small $b>0$ and in any large dimension $n$? Are, conversely, some restrictions on $n$ and/or $a, b$, $\alpha, \beta, u_0$ required? 
%\item c. 
\end{enumerate}
To our knowledge, most of the analyses connected to the aforementioned questions can be found in the literature when the equation for $v$ expressed as (or similarly to) \eqref{KS-NonLocal} is of elliptic type, i.e. for some $\gamma\geq 1$
\begin{equation*}
0=-\Delta v+v +u^\gamma\quad  \textrm{in}\quad \Omega \times (0,\TM). 
\end{equation*}
As a matter of fact, when the equations for the cells and the stimulus are both evolutive, we are only aware of \cite{Tello2021}, where the authors consider, for $\tau=1=m$, $\sigma> 2, \gamma\geq 1$ and $h=h(x,t)\equiv 0$, the initial-boundary value problem associated to this model
\begin{equation}\label{Tello}
\begin{cases}
u_t=\nabla \cdot \left((u+1)^{m-1}\nabla u-\nabla \cdot (\chi u (u+1)^{\sigma-2}\nabla v\right)+f(u) &\quad  \textrm{in}\quad \Omega \times (0,\TM),\\
\tau v_t=\Delta v-v+u^{\gamma}+h &\quad  \textrm{in}\quad \Omega \times (0,\TM).\\
%u_\nu=v_\nu=0 &{\rm on}\ \partial\Omega,\\
%u(x, 0)=u_0(x), v(x,0)=v_0(x), &x \in \Omega,
\end{cases}
\end{equation}
Herein, the nonlocal term is 
\begin{equation}\label{NonLocalTermTello}
f(u):=u\left(a_0-a_1u^\alpha+a_2\int_\Omega u^\alpha dx\right),
\end{equation}
where $\alpha \geq 1$, $a_0,a_1>0$ and $a_2\in \R$; in particular, it is worthwhile mentioning that even though problem \eqref{problem} is the limit case of \eqref{Tello} for $m=1=\gamma$ and $\sigma=2$ (and $h=0$), these models are not directly comparable. In fact, conversely to the mechanism we are dealing with (see again model \eqref{problem}), in \cite{Tello2021} the attractive drift-sensitivity is nonlinear (i.e., $\sigma>2$ in  $-\chi u(u+1)^{\sigma-2} \nabla v$) and, more importantly, the nonlocal term of the reaction in \eqref{NonLocalTermTello} has both an increasing ($a_2>0$) and decreasing ($a_2<0$) effect on the cell density, whereas the dampening counterpart is of polynomial type; this contrasts with \eqref{NonLocalSources}, where the nonlocal term is purely absorbing and the local one productive.  

For model \eqref{Tello} the global-in-time existence of classical solutions and the convergence to the steady state are established in the same \cite{Tello2021}, under suitable regularity assumptions on the initial data and whenever the coefficients of the system satisfy
\begin{equation}\label{ConditionsTello}
\alpha+1>\sigma-1+\gamma\textrm{ and } a_1-a_2|\Omega|>0.
\end{equation}
(Naturally $a_1-a_2|\Omega|>0$ is unnecessary if $a_2\geq 0$.) Additionally, the suppression of some of the conditions in \eqref{ConditionsTello}, might provide (at least from the numerical point of view) some blow-up solution.

As we said above, when the equation for the chemical $v$ is elliptic (biologically this idealizes the situations where chemicals diffuse much faster than cells), some more results are available in the literature. In particular, in \cite{NegreanuTelloNonLocalPara-Ell} the authors analyze, inter alia, problem \eqref{Tello} in the framework of what follows: $\tau=0, \sigma=2$, $m=\gamma=\alpha=1$ and $h=h(x,t)$ is a uniformly bounded function with suitable properties. Similar conclusions as those of the fully parabolic case are derived.  

On the other hand, when the reaction term is taken exactly as in \eqref{NonLocalSources}, these further results dealing with uniform-in-time boundedness of classical solutions emanating from sufficiently regular initial data have been obtained for problem \eqref{Tello}, with $\tau=0$ and $h\equiv 0$: 
%leads to 
\begin{itemize}
\item for the special case where $m=\gamma=a=b=1$ and $\sigma=2$ in \cite{Bian2018}, whenever these assumptions (with $\alpha\geq 1, \beta>1$) $n\geq3,$ $2\leq \alpha <
1+\frac{2\beta}{n}$ or $\frac{n+4}{2}-\beta<\alpha<2$ are complied; 
\item in \cite{Latos2020} for the case  $m=a=b=1$ and $\sigma=2$ $\gamma\geq 1$, $\sigma >2$ tied by $\gamma+\sigma-1\leq \alpha <
1+\frac{2\beta}{n}$ or $\frac{n+4}{2}-\beta<\alpha<\gamma+\sigma-1$;
\item for general choices of the parameters $m>0,\sigma\geq 1,a=b>0$, for $\gamma=1$, under the hypotheses that $\sigma+\frac{n}{2}(\sigma-m)-\beta< \alpha <m+\frac{2}{n}\beta$ or  $\alpha=\sigma+\frac{n}{2}(\sigma-m)-\beta$ together with $b$ large enough (see \cite{TaoFang-NonlinearNonlocal}).
\end{itemize}
For completeness, we add that another indication showing how rich is effectively the study in the framework of models with stationary equations for the stimulus, is given in these papers \cite{BianEtAlWholeSpaceNonLocal,ChenWangDocMat,CarrilloWangAcAppMat,Li-Viglialoro-DIE-Nonlocal}, where nonlocal problems alike those in \eqref{Tello} are studied in the whole space $\R^n$. (In this context, the equation for $v$ is the classical Poisson's equation.)
\subsection{Connection with the Fisher--KPP equation}
In mathematics
\begin{equation}\label{Fisher-KPPeq}
u_t- \Delta u = F(u),
\end{equation} 
is known (in its original one spatial dimensional version) as the Fisher--KPP equation, and it describes a reaction-diffusion phenomenon used to model population growth and wave propagation. (See \cite{FisherPaper,KolmPetroPisk}.) In its more common form $F$, interpretable according to what said above as the rate of growth/death of the population, has this expression ($a,b\geq 0$):
\begin{equation*}
F(u)=a u^\alpha (1 - u)- b u.
\end{equation*}
Apart from the law of the corresponding sources, it appears interesting to discuss the parallelism between equations \eqref{Fisher-KPPeq} and \eqref{KS}: essentially, in the latter the extra transport effect $-\nabla \cdot (u \chi \nabla v)$ appears. In the specific, for $\chi=0$ no convection on the particle density $u$ influences the mechanism, and pure Reaction/$F(u)$-Diffusion/$\Delta u$ models (RDm) are obtained (see \eqref{Fisher-KPPeq}). Oppositely, for $\chi>0$ the population is transported in the habitat toward the direction of $\nabla v$; in this case, equation \eqref{KS} is an example of  Taxis/$\nabla \cdot (u \chi \nabla v)$–Diffusion–Reaction models (TDRm). As a consequence, and at least intuitively,  the sources being equal, TDRm are more inclined to present some instabilities with respect to RDm.

Confining our attention to reactions $F(u)$ of nonlocal type, for a general study on initial-boundary value problems (the majority of them with a homogeneous Dirichlet boundary condition, i.e. $u=0$ on $\partial \Omega$) associated to \eqref{Fisher-KPPeq}, we refer to \cite{quittner-souplet-book,Souplet-BlowUpNonlocal} and references therein. Conversely, for results on more similar contexts to that considered in our analysis, we mention \cite{BianEtAl-FisherKPP}, where  the authors study, among other things, globality and long-time behavior of solutions to a zero-flux nonlocal Fisher–KPP type problem. 
%\begin{equation}\label{problemFisherKKP}
%\begin{cases}
%u_t=\Delta u+ u^\alpha-u^\alpha\int_\Omega u &{\rm in}\ \Omega \times (0, \TM),\\
%%v_t=\Delta v-v+u &{\rm in}\ \Omega \times (0, \TM),\\
%u_\nu=0&{\rm on}\ \partial\Omega \times (0, \TM),\\
%u(x, 0)=u_0(x) &x \in \bar{\Omega},
%\end{cases}
%\end{equation}
%so directly obtainable for $\chi=0$ and $\beta=a=b=1$ from problem \eqref{problem}. It is established that for bounded and smooth domains $\Omega$ of $\R^n$, with $n\in \N$, whenever $u_0=u_0(x)$ is sufficiently regular and such that $\int_\Omega u_0(x)dx<1$, then system \eqref{problemFisherKKP} admits global solutions for $n = 1, 2$ with any $1 \leq \alpha <2$,  or $n \geq 3$ with any $1 \leq  \alpha < 1 + 2/n.$ 
%%%%%%%%%%%%%
\section{Presentation of the main result and organization of the paper}
\subsection{Claim of the main result} 
In this research we intend to improve the degree of knowledge on chemotactic models described by two coupled partial differential equations, and with non-local logistic sources, when both are of parabolic-type. In particular, our overall analysis gives an answer to questions \ref{Question}, in the sense that we establish that \textit{despite any fixed small value of the dampening parameter $b$ and arbitrarily large growth parameter, any initial data $(u_0,v_0)$ (even arbitrarily large) produce  uniform-in-time boundedness of solutions to model \eqref{problem} for both subquadratic and superquadratic growth rate $\alpha$, by properly magnifying the impact associated to the death rate $\beta$.}

Formally, we will prove the following
\begin{theorem}\label{MainTheorem}
Let $\Omega\subset \mathbb{R}^n$, $n\in \N$, be a bounded domain with smooth boundary, $\chi, a, b>0$ and $\alpha, \beta \geq 1$. Additionally, for every  $1<q<\infty$, let $0\leq u_0,v_0\in W^{2,q}(\Omega)$ be given such that $\partial_\nu u_0=\partial_\nu v_0=0$ on $\partial \Omega$. Then, whenever either 
\begin{equation*}
\textnormal{subquadratic growth rate:} \quad 1\leq \alpha <2 \quad \textrm{and} \quad \beta>\frac{n+4}{2}-\alpha,
\end{equation*} 
or
\begin{equation*}
\textnormal{superquadratic growth rate:} \quad \beta>\frac{n}{2} \quad \textrm{and} \quad 2\leq \alpha < 1+ \frac{2\beta}{n},
\end{equation*} 
problem \eqref{problem} admits a unique classical solution, global and uniformly bounded in time, in the sense that 
\begin{equation*}
\begin{cases}
u \in  C^{2,1}(\bar{\Omega}\times(0,\infty))\cap C^{0}(\bar{\Omega}\times[0,\infty)) \cap  L^{\infty}(\bar{\Omega} \times (0, \infty)),&\\
v \in C^{2,1}(\bar{\Omega}\times(0,\infty))\cap C^{0}(\bar{\Omega}\times[0,\infty)) \cap  L^{\infty}_{\textrm{loc}}([0, \infty);W^{1,q}(\Omega)) \cap L^{\infty}(\bar{\Omega} \times (0, \infty)).&
\end{cases}
\end{equation*}
\end{theorem}
\subsection{Structure of the paper} The rest of the paper is structured as follows. First, in $\S$\ref{PreliminarySection}, we collect some necessary and preparatory materials. Then, in $\S$\ref{SectionLocalInTime}, we give some hints on the local-well-posedness to model \eqref{problem}, so obtaining properties of related local solutions $(u,v)$ on $\Omega \times (0,\TM)$; additionally, through the \textit{extensibility criterion} we establish how to ensure globability (i.e., $\TM=\infty$) and boundedness (i.e., $\lVert u(\cdot,t)\rVert_{L^\infty(\Omega)}$ finite on $(0,\infty)$) by using their uniform-in-time $L^k(\Omega)$-boundedness, for $k>1$. Such a bound is derived in $\S$\ref{EstimatesAndProofSection}, and successively used in $\S$\ref{SectionProofTheorem} to prove Theorem \ref{MainTheorem}.
\begin{remark}[On the difficulties of the fully parabolic analysis]
As we will see below, conversely to the parabolic-elliptic case analyzed in \cite[(2.21)]{Bian2018}, in the  fully parabolic case it is no longer possible to use the equation for $v$, so replacing  $\Delta v$ appearing in the testing procedures with $v-u$. This complexity is circumvented by relying on Maximal Sobolev Regularity applied to the equation $v_t=\Delta v-v+u$. 
\end{remark}
\section{Some preliminaries and auxiliary tools}\label{PreliminarySection}
We will make use of this functional relation, obtainable by manipulating the well known Gagliardo--Nirenberg inequality. We underline that for the case $\Omega=\R^n$ 
the proof is given in \cite[Lemma 2]{Bian2016}; we did not find a reference covering bounded domains and henceforth herein we dedicate ourselves to this issue.
\begin{lemma}\label{GagliardoIneqLemma}
Let $\Omega$ be a bounded and smooth domain of $\R^n$, with $n\in \N$ and let, for $n\geq 3$, 
\begin{equation}\label{def_of_p}
p:=\frac{2n}{n-2}. 
\end{equation}
Additionally, let $q, r$ satisfy $1 \le r<q<p$ and $\frac{q}{r}<\frac{2}{r}+1-\frac{2}{p}$. Then for all $\epsilon_1, \epsilon_2>0$ there exists $C_0=C_0(\epsilon_1, \epsilon_2)>0$ such that for all $\varphi \in H^1(\Omega) \cap L^r(\Omega)$, 
\begin{equation}\label{Bian}
\|\varphi\|_{L^q(\Omega)}^q \le C_0 \|\varphi\|_{L^r(\Omega)}^\gamma
+\epsilon_1\|\nabla \varphi\|_{L^2(\Omega)}^2+\epsilon_2\|\varphi\|_{L^2(\Omega)}^2, 
\end{equation}
where 
\begin{equation*}
\lambda:=\frac{\frac{1}{r}-\frac{1}{q}}{\frac{1}{r}-\frac{1}{p}} \in (0,1), \quad \gamma:=\frac{2(1-\lambda)q}{2-\lambda q}.
\end{equation*}
The same conclusion holds for $n\in\{1,2\}$ whenever
$q, r$ fulfill, respectively,  $1 \le r<q$ and $\frac{q}{r}<\frac{2}{r}+2$ and $1 \le r<q$ and $\frac{q}{r}<\frac{2}{r}+1$.
\begin{proof}
Let $n\geq 3$. From the Gagliardo--Nirenberg inequality (\cite[page 126]{Nirenberg1959}) and this algebraic one 
\begin{equation}\label{InequalitiAlgebraic2ToThePower}
(A+B)^l \leq 2^{l-1} (A^l+B^l) \quad \textrm{for all} \quad A,B\geq 0  \; \textrm{ and } \; l\geq 1,
\end{equation}
 for any $q, r>1$ and $s>0$ there is some positive $C_{GN}$ such that  
\begin{equation}\label{GN}
\|\varphi\|_{L^q(\Omega)}^q \leq C_{GN} \|\nabla \varphi\|_{L^2(\Omega)}^{\lambda q} 
\|\varphi\|_{L^r(\Omega)}^{(1-\lambda) q} + C_{GN} \|\varphi\|_{L^s(\Omega)}^{q}, 
\end{equation}
with (recall \eqref{def_of_p})
\begin{equation}\label{lambda}
\lambda=\frac{\frac{1}{r}-\frac{1}{q}}{\frac{1}{r}-\frac{1}{2}+\frac{1}{n}} 
= \frac{\frac{1}{r}-\frac{1}{q}}{\frac{1}{r}-\frac{1}{p}} \in (0,1) \quad \textrm{for all} \quad 
1 \leq r < q < p.
\end{equation}
Now, from the relation  $\frac{q}{r}<\frac{2}{r}+1-\frac{2}{p}$ we have $\frac{\lambda q}{2}<1$, so that the Young inequality applied in \eqref{GN} infers for every $\epsilon_1>0$ some $C_1=C_1(C_{GN}, \epsilon_1) >0$ such that
\begin{equation}\label{GN1}
\|\varphi\|_{L^q(\Omega)}^q \leq \epsilon_1 \|\nabla \varphi\|_{L^2(\Omega)}^2 + C_1 
\|\varphi\|_{L^r(\Omega)}^{\gamma} + C_{GN} \|\varphi\|_{L^s(\Omega)}^{q}, 
\end{equation}
where 
\begin{equation}\label{gamma}
\gamma=\frac{2(1-\lambda)q}{2-\lambda q}.
\end{equation}
On the other hand, for any $q, p >1$, let $s=\frac{2pq}{3p-2}>0$. Subsequently, the H\"{o}lder inequality provides (note that $\frac{2q}{s}=\frac{3p-2}{p}>1$)
\begin{equation*}
C_{GN} \|\varphi\|_{L^s(\Omega)}^{q}= C_{GN} \left(\int_\Omega \varphi^{\frac{s}{q}} \varphi^{s-\frac{s}{q}}\right)^\frac{q}{s} \leq C_{GN} \left(\int_\Omega \varphi^2\right)^{\frac{1}{2}} \left(\int_\Omega \varphi^{\frac{2s(q-1)}{2q-s}}\right)^{\frac{1}{2}(\frac{2q}{s}-1)},
\end{equation*}
and, in turn, Young's inequality gives for any $\epsilon_2>0$, some 
$C_2=C_2(C_{GN}, \epsilon_2) >0$
\begin{equation}\label{H1}
C_{GN} \|\varphi\|_{L^s(\Omega)}^{q} \leq \epsilon_2 \int_\Omega \varphi^2 + C_2 \left(\int_\Omega \varphi^{\frac{2s(q-1)}{2q-s}}\right)^{\frac{2q}{s}-1}.
\end{equation}
The conclusion goes through standard but tedious computations; specifically, by inserting relation \eqref{H1} into estimate \eqref{GN1} and by establishing that for $s$ as above, and $\lambda$ and $\gamma$ as in \eqref{lambda} and \eqref{gamma} respectively, 
$\frac{2s(q-1)}{2q-s}=r$ and $\frac{2q}{s}-1=\frac{\gamma}{r}$, the proof is given with $C_0=C_1+C_2$.

For $n\in\{1,2\}$, the same arguments apply by taking respectively $s=\frac{q}{2}$ and $s=\frac{2q}{3}.$
\end{proof}
\end{lemma}
In the spirit of \cite{IshidaYokotaDCDS-SViaMSR,IshidaYokotaJDE-2012SmallData,SenbaSuzukiAAN2006}, let us recall the following consequence of Maximal Sobolev Regularity results (like \cite{hieber_pruess} or \cite[Thm. 2.3]{giga_sohr}):
\begin{lemma}\label{lem:MaxReg}
 Let $n\in \mathbb{N}$, $\Omega\subset \mathbb{R}^n$ be a bounded and smooth domain and $q\in(1,\infty)$. Then there is $C_{MR}>0$ such that the following holds: Whenever $T\in(0,\infty]$, $I=[0,T)$, $f\in L^q(I;L^q(\Omega))$ and $v_0\in W^{2,q}(\Omega)$ is such that $\partial_{\nu} v_0=0$ on $\partial\Omega$, every solution $v\in W_{loc}^{1,q}(I;L^q(\Omega))\cap L^q_{loc}(I;W^{2,q}(\Omega))$ of
\[
 v_t=\Delta v-v + f\;\; \text{ in }\;\;\Omega\times(0,T);\quad
 \partial_{\nu} v=0\;\; \text{ on }\;\;\partial\Omega \times(0,T); \quad v(\cdot,0)=v_0 \;\; \text{ on }\;\;\Omega
\]
satisfies
\[
 \int_0^t e^s \left(\int_\Omega |\Delta v(\cdot,s)|^q\right)ds \le C_{MR} \left[1+\int_0^t e^s \left(\int_\Omega |f(\cdot,s)|^{q}\right)ds\right] \quad \text{for } 0<t<T.
\]
\begin{proof}
For $(x,t)\in \Omega\times (0,T)$, let us set $z(x,t):=e^{\frac{t}{q}}v(x,t)$.
Then easy computations establish that $z$ solves
\begin{equation*}
\begin{cases}
z_t=\Delta z -\left(1-\frac{1}{q}\right)z +e^\frac{t}{q}f & \textrm{in } \Omega \times (0,T),\\
\partial_\nu z=0 & \textrm{on } \partial \Omega \times (0,T),\\
z(x,0)= v_{0}(x) & x\in  \Omega.
\end{cases}
\end{equation*}
Subsequently, let us apply Maximal Sobolev Regularity (\cite[(3.8)]{hieber_pruess}, \cite[Thm. 2.3]{giga_sohr}) to $A=\Delta-(1-\frac{1}{q})$, $X=L^q(\Omega)$ and $X_1=D(A)=W^{2,q}_\mathcal{\partial_\nu}(\Omega)= \{w\in W^{2,q}(\Omega): \partial_\nu w=0\,\,\textrm{on}\;\partial \Omega\}$, which asserts that with some $c_1>0$ we have for every $t\in(0,T)$ that
\begin{equation*}
\begin{split}
\lVert \Delta z\rVert _{L^q([0,t];L^q(\Omega))} +\lVert z_{t}\rVert_{L^q([0,t];L^{q}(\Omega))}
\leq c_1 \left(\lVert v_{0}\rVert_{1-\frac{1}{q},q}+
\Big(\int_0^t \|e^\frac{s}{q}f(\cdot,s)\|_{L^q(\Omega)}^q\, ds\Big)^\frac{1}{q}\right),
\end{split}
\end{equation*}
where $\lVert \cdot\rVert_{1-\frac{1}{q},q}$ represents the norm in the interpolation space $(X,X_1)_{1-\frac{1}{q},q}$. In turn, we have by using \eqref{InequalitiAlgebraic2ToThePower} that for $C_{MR}=\left(c_1 \max\left\{1,\lVert v_{0}\rVert_{1-\frac{1}{q},q}\right\}\right)^q2^{q-1}$  
\begin{equation}\label{inequalityPrussBIS}
\begin{split}
\int_0^t \Big(\int_\Omega |\Delta z(\cdot,s)|^{q}\Big)\,ds
%&\leq c_1^{q} \left[K+\Big(\int_0^t\|e^\frac{s}{q}f(\cdot,s)\|_{L^{q}(\Omega)}^{q}\, ds\Big)^\frac{1}%{q}\right]^{q}\\ & 
\leq C_{MR}  \left[1+\int_0^t e^s \left(\int_\Omega |f(\cdot,s)|^{q} \right)ds \right]\quad \textrm{for all } t\in(0,T).
\end{split}
\end{equation}
We can finally obtain the claim by re-substituting $z(\cdot,t):=e^{\frac{t}{q}}v(\cdot,t)$ into relation \eqref{inequalityPrussBIS}.
\end{proof}
\end{lemma}
We will also need this comparison argument for Ordinary Differential Equations.
\begin{lemma}\label{LemmaODI-Comparison}
Let\/ $T>0$ and $\phi:(0,T)\times \R^+_0\rightarrow \R$. If $0\leq y\in C^0([0,T))\cap  C^1((0,T))$ is such that 
\begin{equation*}
y'\leq \phi(t,y)\quad \textrm{for all } t \in (0,T), 
\end{equation*}
and there is $y_1>0$ with the property that whenever $y>y_1$ for some $t\in (0,T)$ one has that $\phi(t,y)\leq 0$, then
\begin{equation*}
y\leq \max\{y_1,y(0)\}\quad \textrm{on } (0,T).
\end{equation*}
\begin{proof}
Setting $y_0=y(0)$, let us distinguish the cases $y_0<y_1$ and $y_0\geq y_1$ and let us show that, respectively, the sets
\begin{equation*}
S_{y_1}:=\{t\in (0,T) \mid y(t)>y_1\}\quad \textrm{and} \quad  S_{y_0}:=\{t\in (0,T) \mid y(t)>y_0\}
\end{equation*} 
are empty. In particular, we will establish only that $S_{y_1}=\emptyset$, the reasoning for $S_{y_0}$ being similar. 

By contradiction, if there were some $t_0\in S_{y_1}$ then by the continuity of $y$ and $y_0<y_1$ we could find $I=(\underline{t},\bar{t})$ (with possibly $t_0=\bar{t}$) such that $y_1<y(\underline{t})<y(\bar{t})$, $y_1<y(t)$ on $I$; henceforth, by hypothesis, $\phi(t,y)\leq 0$ for all $t\in I$. At this stage, the Lagrange theorem would provide a proper $\xi \in I$ leading to this inconsistency:
\begin{equation*}
0<\frac{y(\bar{t})-y(\underline{t})}{\bar{t}-\underline{t}}=y'(\xi)\leq \phi(\xi,y)\leq 0.
\end{equation*}
\end{proof}
\end{lemma}
\section{Local solutions and their main properties. A boundedness criterion}\label{SectionLocalInTime}
\begin{lemma}[Local existence and extensibility criterion]\label{localSol}
Let $n\in \mathbb{N}$, $\Omega\subset \mathbb{R}^n$ be a bounded and smooth domain, $\chi, a, b>0$ and $\alpha, \beta \geq 1$. Moreover, for every  $1<q<\infty$, let $u_0,v_0\in W^{2,q}(\Omega)$  satisfy
\begin{equation*}
 \partial_{\nu} u_0 = \partial_{\nu} v_0=0 \mbox{ on } \partial \Omega,  
 \mbox{ and } u_0, v_0\ge 0 \mbox{ on } \bar{\Omega}.
\end{equation*}
Then problem \eqref{problem} has a unique and nonnegative classical solution
\begin{equation*}
\begin{cases}
u \in  C^{2,1}(\bar{\Omega}\times(0,\TM))\cap C^{0}(\bar{\Omega}\times[0,\TM)),&\\
v \in C^{2,1}(\bar{\Omega}\times(0,\TM))\cap C^{0}(\bar{\Omega}\times[0,\TM)) \cap  L^{\infty}_{\textrm{loc}}([0, \TM);W^{1,q}(\Omega)),&
\end{cases}
\end{equation*}
for some maximal $\TM\in(0,\infty]$ which is such that
\begin{equation}\label{extcrit}
 \text{either } \TM=\infty \quad \text{or}\quad \limsup_{t\to \TM} \lVert u(\cdot,t)\rVert_{L^{\infty}(\Omega)} = \infty.
\end{equation}
Additionally, there exists $m_0>0$ such that 
\begin{equation}\label{MassBounded}
\int_\Omega u(x,t)\,dx \leq m_0 \quad \textrm{for all } t \in (0,\TM).
\end{equation}
\begin{proof}
The first part of the proof can be obtained by adapting to the fully parabolic case the reasoning in \cite[Proposition 4]{Bian2018} developed for the simplified parabolic-elliptic scenario. 

As to the boundedness of the mass, we integrate over $\Omega$ the first equation of problem \eqref{problem} so that by H\"{o}lder's inequality, and $\gamma(t):=\int_\Omega u^{\alpha}\geq 0$ on $(0,\TM)$, 
\begin{equation*}%\label{DerivU}
y'(t):=\frac{d}{dt} \int_\Omega u= \int_\Omega u^{\alpha} \left(a-b \int_\Omega u^{\beta}\right)\leq \gamma(t)\left(a-b |\Omega|^{1-\beta} (y(t))^{\beta}\right) 
\quad \textrm{for all } t \in (0, \TM).
\end{equation*}
Now we apply Lemma \ref{LemmaODI-Comparison} with $T=\TM$, $\phi(t,y)=\gamma(t)\left(a-b |\Omega|^{1-\beta} (y(t))^{\beta}\right)$, $y_0=y(0)=\int_\Omega u_0$ and $y_1:=\left(\frac{a}{b |\Omega|^{1-\beta}}\right)^{\frac{1}{\beta}}$, so concluding with $m_0=\max\{y_0, y_1\}$.
\end{proof}
\end{lemma}
Once the classical local well posedness to model \eqref{problem} provided by Lemma \ref{localSol} is ensured (in particular from now on with $(u,v)$ we refer to the local solution defined on $\Omega \times (0,\TM)$), a suitable uniform-in-time boundedness criterion is required. In the specific, the next result based on an iterative method connected to the Moser--Alikakos technique addresses the issue.
\begin{lemma}\label{BoundednessLemma}
Whenever for every $k>1$ there exists $C>0$ such that 
\begin{equation*}
\int_\Omega u^k \leq C \quad \textrm{for all } t \in (0,\TM), 
\end{equation*}
actually $u$ is uniformly bounded, in the sense that $u \in L^{\infty}((0,\infty); L^{\infty}(\Omega))$. Automatically, $v$ is also uniformly bounded.
\begin{proof}
From the first equation of problem \eqref{problem} and the nonnegativity of $u$, we have that $u$ itself is such that $u_t \leq \Delta u -\chi \nabla \cdot (u \nabla v) + a u^{\alpha}$. In particular, $u$ solves \cite[(A.1)]{TaoWinkParaPara} with $D(x,t,u)=1$, $f(x,t)=-\chi u(x,t) \nabla v(x,t)$ and $g(x,t)= a u^{\alpha}(x,t)$. In these positions, since from our hypotheses $u \in L^{\infty}((0,\TM); L^k(\Omega))$ for all $k>1$ (and in particular for $k$ arbitrarily large),  $g$ belong to $L^{\infty}((0,\TM); L^k(\Omega))$ and from parabolic regularity results (\cite[IV. 5.3]{LSUBookInequality}) we have that also $\nabla v \in L^{\infty}((0,\TM); L^k(\Omega))\in L^{\infty}((0,\TM); L^k(\Omega))$.
As a by-product, $f$ and, and
\cite[Lemma A.1]{TaoWinkParaPara} ensures $u \in L^{\infty}((0,\TM); L^\infty(\Omega))$.
Finally, the extensibility criterion \eqref{extcrit} entails $\TM=\infty$ and we conclude. (The boundedness of $v$ follows from $u\in L^{\infty}((0,\infty); L^k(\Omega))$ for arbitrarily large  $k>1$ and, again, parabolic regularity results and Sobolev embeddings.)
\end{proof} 
\end{lemma}
\section{A priori estimates}\label{EstimatesAndProofSection}
Since the uniform-in-time boundedness of $u$ is implied whenever $u\in L^\infty((0,\TM);L^k(\Omega))$ for some $k>1$, here under we dedicate to the derivation of some \textit{a priori} integral estimates. 

(\textit{In the sequel we will tacitly assume that all the constants $c_i$ appearing below, $i=1,2,\ldots$ are positive.}) 
\begin{lemma}\label{Est_u_Lap_v}
For all $k>1$, $\chi>0$ there exist $\const{a}, \const{b}$ such that 
whenever $\alpha>1$
\begin{equation}\label{AlphaMAg1}
(k-1) \chi \int_\Omega u^k\Delta v \le 
\int_\Omega u^{k+\alpha-1}+\const{a}\int_\Omega |\Delta v|^{\frac{k+\alpha-1}{\alpha-1}}\quad \textrm{for all } t \in (0, \TM),
\end{equation}
wile if $\alpha\geq 1$.
\begin{equation}\label{AlphaMAgUgu1}
(k-1) \chi \int_\Omega u^k\Delta v \leq
\int_\Omega u^{k+1}+\const{b}\int_\Omega |\Delta v|^{k+1} \quad \textrm{for all } t \in (0, \TM).
\end{equation}
\begin{proof}
The Young inequality directly provides the claim. 
\end{proof}
\end{lemma}
Let us now distinguish the analysis of the subquadratic case from the superquadratic one, exactly starting from this last situation.
%\begin{lemma}\label{para_reg}
%For $f \in L^p((0, T); L^q(\Omega))$, $1<p, q<\infty$, $T<\infty$ and $\varphi_0 \in W^{2,q}(\Omega)$, such that
%\begin{equation*}
%(\varphi_0)_\nu=0\quad {\rm on}\ \partial \Omega,
%\end{equation*}
%the solution $\varphi$ of the system 
%\begin{equation*}
%\begin{cases}
%\varphi_t=\Delta \varphi-\varphi+f & {\rm in}\ \Omega \times (0,T),\\
%\varphi_\nu=0 & {\rm on}\ \partial \Omega,\\
%\varphi(x,0)=\varphi_0(x) & x \in \Omega,
%\end{cases}
%\end{equation*}
%satisfies $\varphi \in W^{1,p}((0,T); L^q(\Omega)) \cap L^p((0,T); W^{2,q}(\Omega))$. Moreover, for $p=q$, the solution $\varphi$ fulfills 
%\begin{equation*}
%\int_0^T e^t \int_\Omega (|D_x^2\varphi|^q+|D_x\varphi|^q+|\varphi|^q)\,dxdt \le \const{A}\left(\int_0^T e^t\int_\Omega f^q\,dxdt +\|\varphi_0\|_{W^{2,q}(\Omega)}^q\right).
%\end{equation*}
%\end{lemma}
\subsection{The superquadratic growth: $\beta>\frac{n}{2}$ and $2 \le \alpha<1+\frac{2\beta}{n}$}
\begin{lemma}\label{CaseAlphaMaggiore2Lemma}
Assume that $\alpha, \beta \ge 1$ satisfy that
\begin{equation}\label{alpha_large}
\beta > \frac{n}{2} \quad \textrm{and} \quad 2 \le \alpha<1+\frac{2\beta}{n}.
\end{equation}
Then there exist $k_0\geq1, L_0>0$ such that for all $k>k_0$, 
\begin{equation*}
\int_\Omega u^k \le L_0 \quad  \mbox{for\ all}\ t \in (0, \TM). 
\end{equation*}
\begin{proof}
Let us start fixing $k_0=1$, and when necessary we will enlarge this initial value. For all $k>k_0$, we have from the first equation in \eqref{problem} and integration by parts that 
\begin{equation}\label{uk_estimate1}
\begin{split}
\frac{d}{dt}\int_\Omega u^k &=k\int_\Omega u^{k-1}\Delta u-k\chi \int_\Omega u^{k-1}\nabla \cdot (u\nabla v) +ka\int_\Omega u^{k+\alpha-1}\\
&\quad -kb\left(\int_\Omega u^{k+\alpha-1}\right)\left(\int_\Omega u^\beta\right)\\
&=-k(k-1)\int_\Omega u^{k-2}|\nabla u|^2+k(k-1)\chi\int_\Omega u^{k-1}\nabla u\cdot \nabla v +ka\int_\Omega u^{k+\alpha-1}\\
&\quad -kb\left(\int_\Omega u^{k+\alpha-1}\right)\left(\int_\Omega u^\beta\right)\\
&=-\frac{4(k-1)}{k}\int_\Omega |\nabla u^{\frac{k}{2}}|^2-(k-1)\chi\int_\Omega u^k\Delta v +ka\int_\Omega u^{k+\alpha-1}\\
&\quad -kb\left(\int_\Omega u^{k+\alpha-1}\right)\left(\int_\Omega u^\beta\right) 
\quad  \mbox{on}\ (0, \TM). 
\end{split}
\end{equation}
Here, from bound \eqref{AlphaMAg1} in Lemma~\ref{Est_u_Lap_v} we have that 
\begin{equation}\label{uLap_v_estimate}
-(k-1)\chi\int_\Omega u^k\Delta v \le \int_\Omega u^{k+\alpha-1} + \const{a} \int_\Omega |\Delta v|^{\frac{k+\alpha-1}{\alpha-1}} \quad  \mbox{for\ all}\ t \in (0, \TM).
\end{equation}
A combination of relations \eqref{uk_estimate1} and \eqref{uLap_v_estimate} implies that
for all $t \in (0, \TM)$
\begin{equation}\label{diff_u^k}
\frac{d}{dt}\int_\Omega u^k+kb\left(\int_\Omega u^{k+\alpha-1}\right)\left(\int_\Omega u^\beta\right) \le -\frac{4(k-1)}{k}\int_\Omega |\nabla u^{\frac{k}{2}}|^2
+\const{g} \int_\Omega u^{k+\alpha-1}+\const{a}\int_\Omega|\Delta v|^{\frac{k+\alpha-1}{\alpha-1}}.
\end{equation}
We now estimate the second integral on the right-hand side of \eqref{diff_u^k}.
From the identity $\int_\Omega u^{k+\alpha-1}=\|u^{\frac{k}{2}}\|_{L^{\frac{2(k+\alpha-1)}{k}}(\Omega)}^{\frac{2(k+\alpha-1)}{k}}$, our aim is exploiting Lemma \ref{GagliardoIneqLemma} with $\varphi:=u^{\frac{k}{2}}$ and proper $q$ and $r$. In the specific, for $n\geq 3$ (at the end of this proof we will discuss the cases $n=1$ and $n=2$) in order to make meaningful the forthcoming computations, let us take $k_0=\max\{\beta-\alpha+1,1\}$. From the definition of $k_0$ and condition \eqref{alpha_large}, for any $k>k_0$ it is possible to set 
\begin{equation}\label{kPrimo}
k':=\frac{k+\alpha+\beta-1}{2},
\end{equation}
which satisfies 
\begin{equation}\label{k'condi}
\max\left\{\beta,\ \frac{k}{2},\ \frac{p(\alpha-1)}{p-2}\right\}<k'<k+\alpha-1.
\end{equation}
In this way, for 
\begin{equation*}
q:=\frac{2(k+\alpha-1)}{k},\ r:=\frac{2k'}{k}
\end{equation*}
a number of calculations yield $1\le r<q<p$ and $\frac{q}{r}<\frac{2}{r}+1-\frac{2}{p}$. Therefore we infer from \eqref{Bian} that for all $\bar{c}>0$
\begin{equation}\label{u^k+alpha-1}
\bar{c} \int_\Omega u^{k+\alpha-1}= \bar{c}\|u^{\frac{k}{2}}\|_{L^{\frac{2(k+\alpha-1)}{k}}(\Omega)}^{\frac{2(k+\alpha-1)}{k}}
\le \frac{2(k-1)}{k}\int_\Omega |\nabla u^{\frac{k}{2}}|^2+\int_\Omega u^k+\const{i}\left(\int_\Omega u^{k'}\right)^{\frac{\gamma}{r}}\quad  \mbox{for\ all}\ t \in (0, \TM). 
\end{equation}
Here, the interpolation inequality (see \cite[page 93]{Brezis}) yields for all $t \in (0, \TM)$, 
\begin{equation}\label{u^kp}
\begin{split}
\left(\int_\Omega u^{k'}\right)^{\frac{\gamma}{r}}=\|u\|_{L^{k'}(\Omega)}^{b_1}
&\le \|u\|_{L^\beta(\Omega)}^{a_1b_1} \|u\|_{L^{k+\alpha-1}(\Omega)}^{(1-a_1)b_1}\\
&=\left(\|u\|_{L^\beta(\Omega)}^\beta \|u\|_{L^{k+\alpha-1}(\Omega)}^{k+\alpha-1}\right)^{\frac{(1-a_1)b_1}{k+\alpha-1}}\|u\|_{L^\beta(\Omega)}
^{\left[a_1-\frac{\beta(1-a_1)}{k+\alpha-1}\right]b_1},
\end{split}
\end{equation}
where 
\begin{equation}\label{DefinitionB1A1}
b_1=b_1(q):=\frac{k'\gamma(q)}{r}=\frac{k'\gamma}{r}, \quad a_1:=\frac{\frac{1}{k'}-\frac{1}{k+\alpha-1}}{\frac{1}{\beta}-\frac{1}{k+\alpha-1}}\in (0,1).
\end{equation}
We note that recalling the expression of $k'$ in \eqref{kPrimo} and the range of $\alpha$ in \eqref{alpha_large}, some computations provide 
\begin{equation*}%\label{pow_zero}
\left[a_1-\frac{\beta(1-a_1)}{k+\alpha-1}\right]b_1=0 \quad \textrm{and} \quad 
\frac{(1-a_1)b_1}{k+\alpha-1}<1. 
\end{equation*}
As a consequence, we can invoke Young's inequality so that  
relation \eqref{u^kp} reads 
\begin{equation*}
\const{i}\left(\int_\Omega u^{k'}\right)^{\frac{\gamma}{r}} \le \const{i}\left(\|u\|_{L^\beta(\Omega)}^\beta \|u\|_{L^{k+\alpha-1}(\Omega)}^{k+\alpha-1}\right)
^{\frac{(1-a_1)b_1}{k+\alpha-1}}
\le kb\left(\int_\Omega u^{k+\alpha-1}\right)\left(\int_\Omega u^\beta\right)+\const{j} \quad  \mbox{for\ all}\ t \in (0, \TM),
\end{equation*}
which in conjunction with \eqref{u^k+alpha-1} implies for all $t \in (0, \TM)$, 
\begin{equation}\label{u^k+alpha-1_2}
\const{g} \int_\Omega u^{k+\alpha-1} \le \frac{2(k-1)}{k}\int_\Omega |\nabla u^{\frac{k}{2}}|^2+\int_\Omega u^k
+kb\left(\int_\Omega u^{k+\alpha-1}\right)\left(\int_\Omega u^\beta\right)+\const{j}.
\end{equation}
Now we focus on the second integral at the right-hand side: the Gagliardo--Nirenberg inequality and \eqref{MassBounded} produce for
\begin{equation*}
\theta_1:=\frac{\frac{k}{2}-\frac{1}{2}}{\frac{k}{2}+\frac{1}{n}-\frac{1}{2}} \in (0,1),
\end{equation*}
this bound on $(0,\TM)$:
\begin{equation*}
\int_\Omega u^k=\|u^\frac{k}{2}\|_{L^2(\Omega)}^2
\le \const{l} \|\nabla u^{\frac{k}{2}}\|_{L^2(\Omega)}^{2\theta_1}\|u^\frac{k}{2}\|_{L^\frac{2}{k}(\Omega)}^{2(1-\theta_1)}
+\const{l}\|u^\frac{k}{2}\|_{L^\frac{2}{k}(\Omega)}^2
\le \const{m}\left(\int_\Omega |\nabla u^{\frac{k}{2}}|^2\right)^{\theta_1}+\const{m}.
\end{equation*}
In turn, we have from the Young inequality that for all $\hat{c}>0$
\begin{equation}\label{Young}
\begin{split}
\hat{c} \int_\Omega u^k &\le \frac{2(k-1)}{k}\int_\Omega |\nabla u^{\frac{k}{2}}|^2+\const{n} \quad  \mbox{for\ all}\ t \in (0, \TM). 
\end{split}
\end{equation}
Coming back to \eqref{diff_u^k}, in order to estimate the term $\const{a}\int_\Omega|\Delta v|^{\frac{k+\alpha-1}{\alpha-1}}$, let us exploit Lemma \ref{lem:MaxReg} with $q=\frac{k+\alpha-1}{\alpha-1}$: we have 
\begin{equation}\label{EstimateLaplV}
\const{a} \int_0^t e^s \left(\int_\Omega |\Delta v(\cdot,s)|^{\frac{k+\alpha-1}{\alpha-1}}
\right)ds \le \const{a} C_{MR} \left[1+\int_0^t e^s \left(\int_\Omega u(\cdot,s)^{\frac{k+\alpha-1}{\alpha-1}}\right)ds\right] \quad \text{for all } t\in(0,\TM).
\end{equation}
Since from the condition $\alpha \ge 2$ we have that $\frac{k+\alpha-1}{\alpha-1} \le k+\alpha-1$, the Young inequality leads to
\begin{equation}\label{EstU}
\const{a} C_{MR} \int_\Omega u^{\frac{k+\alpha-1}{\alpha-1}} \le \const{a} C_{MR}\int_\Omega u^{k+\alpha-1}+\const{p} \quad  \mbox{for\ all}\ t \in (0, \TM).
\end{equation}
(Naturally for the limit case $\alpha=2$, the constant $\const{p}$ can be taken equal to $0$.)
We now add to both sides of \eqref{diff_u^k} the term $\int_\Omega u^k$ and then 
we multiply by $e^t$. Since $e^t \frac{d}{dt}\int_\Omega u^k + e^t \int_\Omega u^k =
\frac{d}{dt}\left(e^t \int_\Omega u^k\right)$, an integration over $(0,t)$ provides for all $t \in (0,\TM)$
\begin{equation}\label{diff1_u^k}
\begin{split}
&e^t \int_\Omega u^k - \int_\Omega u^k_0 +kb \int_0^t e^s \left(\int_\Omega u^{k+\alpha-1}\right)\left(\int_\Omega u^\beta\right)\,ds
\le -\frac{4(k-1)}{k} \int_0^t e^s  \left(\int_\Omega |\nabla u^{\frac{k}{2}}|^2\right)\,ds\\
&+\int_0^t e^s \left(\int_\Omega u^k\right)\,ds
+\const{g} \int_0^t e^s \left(\int_\Omega u^{k+\alpha-1}\right)\,ds +\const{a} \int_0^t e^s \left(\int_\Omega|\Delta v|^{\frac{k+\alpha-1}{\alpha-1}}\right)\,ds.
\end{split}
\end{equation}
By inserting estimate \eqref{EstimateLaplV} into \eqref{diff1_u^k} and taking into account bounds \eqref{EstU}, \eqref{u^k+alpha-1_2} and \eqref{Young}, we arrive at
%\begin{equation}\label{diff2_u^k}
%\begin{split}
%&\int_0^t e^s \frac{d}{dt}\int_\Omega u^k+ (kb-2\epsilon_3) \int_0^t e^s \left(\int_\Omega u^{k+\alpha-1}\right)\left(\int_\Omega u^\beta\right)\,ds\\
%&\le \left(-\frac{4(k-1)}{k}+2\epsilon_1+2\epsilon_4\right) \int_0^t e^s  \left(\int_\Omega |\nabla u^{\frac{k}{2}}|^2\right)\,ds +\const{gs}e^t + \const{gs1}  \quad \textrm{on } (0, \TM). 
%\end{split}
%\end{equation}
\begin{equation*}%\label{diff2_u^k}
\begin{split}
e^t \int_\Omega u^k \leq \int_\Omega u^k_0 +\const{gs}e^t + \const{gs1}  \quad \textrm{on } (0, \TM),
\end{split}
\end{equation*}
which implies 
\begin{equation*}
\int_\Omega u^k \le L_0\quad  \mbox{for\ all}\ t \in (0, \TM)
\end{equation*}
with $L_0:=\const{gs2}+\int_\Omega u_0^k$, so the claim is proved. 

For $n\in\{1,2\}$ the arguments are similar once relation \eqref{k'condi} is, respectively, replaced by 
$$
\max\left\{\beta,\frac{k}{2},\frac{\alpha-1}{2}\right\}<k'<k+\alpha -1\quad \textrm{and} \quad \max\left\{\beta,\frac{k}{2},\alpha-1\right\}<k'<k+\alpha -1. 
$$ 
\end{proof}
\end{lemma}
\subsection{The subquadratic growth: $1 \le \alpha <2$ and $\beta>\frac{n+4}{2}-\alpha$}
\begin{lemma}\label{CaseAlphaMinore2Lemma}
Assume that $\alpha, \beta \ge 1$ satisfy 
\begin{equation}\label{beta_large}
1 \le \alpha <2 \quad \mbox{and}\quad \beta>\frac{n+4}{2}-\alpha. 
\end{equation}
Then there exist $k_1\geq 1, L_1>0$ such that for all $k>k_1$, 
\begin{equation*}
\int_\Omega u^k \le L_1 \quad  \mbox{for\ all}\ t \in (0, \TM). 
\end{equation*}
\begin{proof}
Let us consider $k_1=1$; as done before, we will enlarge this initial value when necessary. By following the same argument of Lemma \ref{CaseAlphaMaggiore2Lemma}
for all $k>k_1$, we arrive for all $t \in (0, \TM)$ at
\begin{equation}\label{uk_estimate2}
\frac{d}{dt}\int_\Omega u^k =-\frac{4(k-1)}{k}\int_\Omega |\nabla u^{\frac{k}{2}}|^2-(k-1)\chi\int_\Omega u^k\Delta v +ka\int_\Omega u^{k+\alpha-1}-kb\left(\int_\Omega u^{k+\alpha-1}\right)\left(\int_\Omega u^\beta\right).
\end{equation}
Since $\alpha\geq 1$, an application of relation \eqref{AlphaMAgUgu1} of Lemma~\ref{Est_u_Lap_v} to the second integral at the right-hand side of \eqref{uk_estimate2} gives 
\begin{equation}\label{Lem14}
-(k-1)\chi\int_\Omega u^k\Delta v \le \int_\Omega u^{k+1}+\const{b}\int_\Omega |\Delta v|^{k+1} \quad  \mbox{for\ all}\ t \in (0, \TM),
\end{equation}
whereas from the condition $\alpha<2$, the Young inequality leads to
\begin{equation}\label{Young14}
ka\int_\Omega u^{k+\alpha-1} \le \int_\Omega u^{k+1}+\const{w} \quad  \mbox{for\ all}\ t \in (0, \TM).
\end{equation}
Combining estimates \eqref{Lem14} and \eqref{Young14} with bound \eqref{uk_estimate2}, we have for all $t \in (0, \TM)$, 
\begin{equation}\label{uk_estimate3}
\frac{d}{dt}\int_\Omega u^k +kb\left(\int_\Omega u^{k+\alpha-1}\right)\left(\int_\Omega u^\beta\right) \leq -\frac{4(k-1)}{k}\int_\Omega |\nabla u^{\frac{k}{2}}|^2 + 2 \int_\Omega u^{k+1}+\const{b}\int_\Omega |\Delta v|^{k+1}+\const{w}.
\end{equation}
Now let us focus on the second integral on the right-hand side of \eqref{uk_estimate3}. Since $\int_\Omega u^{k+1}=\|u^{\frac{k}{2}}\|_{L^{\frac{2(k+1)}{k}}(\Omega)}^{\frac{2(k+1)}{k}}$, we can apply Lemma \ref{GagliardoIneqLemma} with $\varphi:=u^{\frac{k}{2}}$ and suitable $q$ and $r$. In the specific, for any 
\begin{equation*}%\label{kcondi2}
k>k_1:=\max\left\{1,1-\alpha+\beta\right\},
\end{equation*}
by posing 
\begin{equation*}
k':=\frac{k+\alpha+\beta-1}{2}, 
\end{equation*}
it is possible to check that 
\begin{equation}\label{k'condi2}
\max\left\{\beta, \frac{k}{2}, \frac{p}{p-2}\right\}<k'<k+\alpha-1.
\end{equation}
In this way, and for $n\geq 3$, letting
\begin{equation*}
q:=\frac{2(k+1)}{k},\ r:=\frac{2k'}{k}
\end{equation*}
we can establish that $1 \le r<q<p$ and $\frac{q}{r}<\frac{2}{r}+1-\frac{2}{p}$. Consequently,  we deduce from \eqref{Bian} that for all $\tilde{c}>0$
\begin{equation}\label{u^k+1}
\tilde{c}\|u^{\frac{k}{2}}\|_{L^{\frac{2(k+1)}{k}}(\Omega)}^{\frac{2(k+1)}{k}} \le \frac{2(k-1)}{k}\int_\Omega |\nabla u^{\frac{k}{2}}|^2+\int_\Omega u^k+\const{y}\left(\int_\Omega u^{k'}\right)^{\frac{\gamma}{r}}\quad  \mbox{for\ all}\ t \in (0, \TM).
\end{equation}
Now an application of the interpolation inequality yields for all $t \in (0, \TM)$, 
\begin{equation*}%\label{u^kp2}
\begin{split}
\left(\int_\Omega u^{k'}\right)^{\frac{\gamma}{r}}=\|u\|_{L^{k'}(\Omega)}^{b_2} 
&\le \|u\|_{L^\beta(\Omega)}^{a_2b_2} \|u\|_{L^{k+\alpha-1}(\Omega)}^{(1-a_2)b_2}\\
&=\left(\|u\|_{L^\beta(\Omega)}^\beta \|u\|_{L^{k+\alpha-1}(\Omega)}^{k+\alpha-1}\right)^{\frac{a_2b_2}{\beta}}\|u\|_{L^{k+\alpha-1}(\Omega)}
^{\left[1-a_2-\frac{a_2(k+\alpha-1)}{\beta}\right]b_2},
\end{split}
\end{equation*}
where 
\begin{equation*}
b_2=b_2(q):=\frac{k'\gamma(q)}{r}=\frac{k'\gamma}{r}, \quad 
a_2:=\frac{\frac{1}{k'}-\frac{1}{k+\alpha-1}}{\frac{1}{\beta}-\frac{1}{k+\alpha-1}} \in (0,1).
\end{equation*}
(A comparison between the couple $(a_2,b_2)$ above and $(a_1,b_1)$ in \eqref{DefinitionB1A1} shows that $a_1=a_2$, whereas $b_i$, $i=1,2$ depends on $q$.)
From straightforward calculations and the condition \eqref{beta_large}, we observe 
that 
\begin{equation*}%\label{pow_zero2}
\left[1-a_2-\frac{a_2(k+\alpha-1)}{\beta}\right]b_2=0 \quad \textrm{and} \quad \frac{a_2b_2}{\beta}<1. 
\end{equation*}
Subsequently, we can exploit the Young inequality entailing 
\begin{equation*}
\const{y}\left(\int_\Omega u^{k'}\right)^{\frac{\gamma}{r}}\le\const{y}\left(\|u\|_{L^\beta(\Omega)}^\beta\|u\|_{L^{k+\alpha-1}(\Omega)}^{k+\alpha-1}\right)^{\frac{a_2b_2}{\beta}}
\le kb\left(\int_\Omega u^{k+\alpha-1}\right)\left(\int_\Omega u^\beta\right)+\const{z} \quad  \mbox{on }\;(0, \TM).
\end{equation*}
This in conjunction with \eqref{u^k+1} implies that 
$t \in (0, \TM)$, 
\begin{equation}\label{u^k+1_2}
\tilde{c} \int_\Omega u^{k+1}\le \frac{2(k-1)}{k} \int_\Omega |\nabla u^{\frac{k}{2}}|^2+
 \int_\Omega u^k +
k b \left(\int_\Omega u^{k+\alpha-1}\right)\left(\int_\Omega u^\beta\right)
+\const{z}. 
\end{equation}
As to the term $\int_\Omega |\Delta v|^{k+1}$ in expression \eqref{uk_estimate3},
by exploiting Lemma \ref{lem:MaxReg} with $q=k+1$, we obtain 
\begin{equation}\label{Estimate1LaplV}
\const{b} \int_0^t e^s \left(\int_\Omega |\Delta v(\cdot,s)|^{k+1}
\right)ds \le \const{b} C_{MR} \left[1+\int_0^t e^s \left(\int_\Omega u(\cdot,s)^{k+1} \right)ds\right] \quad \text{for all } t\in(0,\TM).
\end{equation}
On the other hand, by adding $\int_\Omega u^k$ at both sides of estimate \eqref{uk_estimate3}, by multiplying what obtained by $e^t$, a subsequent integration over $(0,t)$ yields 
\begin{equation}\label{diff_u^k 143}
\begin{split}
&e^t\int_\Omega u^k-\int_\Omega u_0^k +kb \int_0^t e^s \left(\int_\Omega u^{k+\alpha-1}\right)\left(\int_\Omega u^\beta\right)\,ds\\
&\leq -\frac{4(k-1)}{k} \int_0^t e^s \left(\int_\Omega |\nabla u^{\frac{k}{2}}|^2\right)\,ds 
+2 \int_0^t e^s \left(\int_\Omega u^{k+1}\right)\,ds +  \int_0^t e^s \left(\int_\Omega u^k\right)\,ds\\
&\quad +\const{b}\int_0^t e^s \left(\int_\Omega|\Delta v|^{k+1}\right)\,ds +\const{aa}e^t \quad  \mbox{for\ all}\ t \in (0, \TM).
\end{split}
\end{equation}
By rearranging bound \eqref{diff_u^k 143} by virtue of estimates \eqref{Estimate1LaplV}, \eqref{u^k+1_2} and \eqref{Young}, it is provided
\begin{equation*}
e^t \int_\Omega u^k \leq \int_\Omega u_0^k + \const{ad} e^t + \const{ad1} \quad \textrm{on }\, (0,\TM),
\end{equation*}
which gives  
\begin{equation*}
\int_\Omega u^k \le L_1\quad  \mbox{for\ all}\ t \in (0, \TM)
\end{equation*}
with $L_1:=\const{ad2}+ \int_\Omega u_0^k$, so proving the claim.  

To establish the claim for $n\in\{1,2\}$, relation \eqref{k'condi2} has to be taken as
$$
\max\left\{\beta,\frac{k}{2}\right\}<k'<k+\alpha -1. 
$$ 
\end{proof}
\end{lemma}

\section{Proof of Theorem \ref{MainTheorem}}\label{SectionProofTheorem}
We apply Lemma \ref{CaseAlphaMaggiore2Lemma} and Lemma \ref{BoundednessLemma}, and Lemma \ref{CaseAlphaMinore2Lemma} and Lemma \ref{BoundednessLemma} to give the proof for the subquadratic and superquadratic case, respectively. 
\qed

\subsubsection*{\bf\textit{\quad Acknowledgments}}
SF and GV are members of the Gruppo Nazionale per l'Analisi Matematica, la Probabilit\`a e le loro Applicazioni (GNAMPA) of the Istituto Nazionale di Alta Matematica (INdAM), and are partially supported by the research project {\em  Analysis of PDEs in connection with real phenomena} (2021, Grant Number: F73C22001130007), funded by  \href{https://www.fondazionedisardegna.it/}{Fondazione di Sardegna}. 
GV is also supported by MIUR (Italian Ministry of Education, University and Research) Prin 2022 \emph{Nonlinear differential problems with applications to real phenomena} (Grant Number: 2022ZXZTN2), and 
acknowledges financial support under the National Recovery and Resilience Plan (NRRP), Mission 4 Component 2 Investment 1.5 - Call for tender No.3277 published on December 30, 2021 by the Italian Ministry of University and Research (MUR) funded by the European Union -- NextGenerationEU. Project Code ECS0000038 -- Project Title eINS Ecosystem of Innovation for Next Generation Sardinia -- CUP F53C22000430001- Grant Assignment Decree No. 1056 adopted on June 23, 2022 by the Italian Ministry of University and Research (MUR).

%\bibliography{Bibliography}{}
%\bibliographystyle{abbrv}
%\bibliographystyle{alpha}
\end{document}